\documentclass[12pt,reqno]{amsart}
\usepackage{color}
\usepackage{epsfig}
\usepackage{amsmath}
\usepackage{mathtools}
\usepackage{amssymb}
\usepackage{amsthm}
\usepackage{amscd}
\usepackage{graphicx}
\swapnumbers                  
\usepackage[latin1]{inputenc}
\usepackage[T1]{fontenc}
\usepackage[french, english]{babel}
\usepackage[totalwidth=400pt, totalheight=620pt]{geometry}

\newcommand{\C}{{\mathbb{C}}}

\newcommand{\N}{{\mathbb{N}}}

\newcommand{\ccup}{\cup\hspace{-6pt}{^\uparrow}}

\newcommand{\dref}[1]{Definition~\ref{#1}} 
\newcommand{\eref}[1]{~(\ref{#1})}

\newcommand{\lref}[1]{Lemma~\ref{#1}}
\newcommand{\pref}[1]{Proposition~\ref{#1}}
\newcommand{\rref}[1]{Remark~\ref{#1}}
\newcommand{\tref}[1]{Theorem~\ref{#1}}
\newcommand{\sref}[1]{Section~\ref{#1}}

\def\preuve{\begin{proof}}

\newtheorem{defi}[subsubsection]{Definition}
\newtheorem{Def}[subsection]{Definition}

\newtheorem{defini}[subsection]{Definitions}

\newtheorem{lemm}[subsubsection]{Lemma}
\newtheorem{prop}[subsubsection]{Proposition}
\newtheorem{pro}[subsection]{Proposition}
\newtheorem{rem}[subsubsection]{Remark}

\newtheorem{theo}[subsubsection]{Theorem}

\newtheorem{exem}[subsubsection]{Example}

  \newenvironment{demo}{\noindent {\it Proof:}
      \begin{quotation}\noindent}{\end{quotation}\hfill$\square $}

\let\<\langle
\let\>\rangle
\let\uml\"

\title[\MakeLowercase{m}-accretive Laplacian on a non symmetric graph]{\MakeLowercase{m}-accretive Laplacian on \\a non symmetric graph}

\author[Colette Ann\'e]{Colette Ann\'e}
\address{Universit\'e de Nantes, Laboratoire de Math\'ematique Jean Leray, CNRS, Facult\'e des Sciences, BP 92208, 44322 Nantes, (France).}

\email{colette.anne@univ-nantes.fr}

\author{Marwa Balti}

\address{Universit\'e de Monastir, (LR/18ES15) (Tunisie).\\
Universit\'e de Nantes, Laboratoire de Math\'ematique Jean Leray, CNRS, Facult\'e des Sciences, BP 92208, 44322 Nantes, (France).
}
\email{balti-marwa@hotmail.fr}

\author{Nabila Torki-Hamza}
\address{Universit\'e de Monastir, LR/18ES15 \& Institut Sup\'erieur d'Informatique de Mahdia (ISIMa) B.P 49, Campus Universitaire de Mahdia; 5111-Mahdia (Tunisie).}
\email{natorki@gmail.com}

\date{\today, \emph{File:} \texttt{\jobname.tex}}


\begin{document}
\begin{abstract}
  We consider a non self-adjoint Laplacian on a directed graph with non symmetric weights on edges. We give a
  criterion for the m-accretiveness and the m-sectoriality of this Laplacian. Our results are based on
  a comparison of this operator with its symmetric part for which we can apply different results concerning
  essential self-adjointness of a symmetric Laplace operator on an infinite graph. This gives results on the
  heat operator related to our non-symmetric Laplacian.
\end{abstract}
\maketitle
\tableofcontents

\section*{Introduction}
Some properties for linear unbounded non self-adjoint operators, such as accretiveness, maximal accretiveness (or
m-accretiveness) and m-sectoriality are very important for physical and technical problems.  They are subjects of
special attention in view of later applications to analytic and asymptotic perturbation
theory \cite{Kts}, \cite{Ev}, \cite{Helf}, \cite{Ou}, \cite{Kh}, \cite{Mil}. The main importance of accretive
operators is their appearance in the Hille-Yosida and Lumer-Phillips Theorems: an operator $A$ is maximally
accretive if and only if $-A$ is the generator of a contraction semigroup. Also, we shall focus on
m-sectorial operators:
their spectrum lies in a sector and their resolvent satisfies a certain estimate.
The opposite of generators of bounded holomorphic semigroups holomorphic on a sector are m-sectorial
operators.\\

We consider a directed infinite graph and we investigate the associated non symmetric Laplacian $\Delta$ under
a Kirchhoff's Assumption \cite{Mar}. This class of operators can be considered as a generalization of
lower semibounded or positive symmetric operators \cite{To10}, \cite{kel}. The purpose of the present paper is
to give a criterion for the m-accretiveness and the m-sectoriality of the discrete Laplacian.

After the preliminaries in Section 1, we have four sections. In Section 2, we give a general condition on
the graph using the notion of $\chi$-completeness introduced in \cite{AT} to establish the m-accretiveness of
the non symmetric Laplacian $\Delta$. Section 3 is devoted to the study of the relations between our non
symmetric Laplacian and its symmetrized part. Section 4 deals with the m-sectoriality of $\Delta$ as a generalisation
of \cite{A}. Section 5 presents properties induced by the  m-accretivity of our operator.
\section{Preliminaries}\label{In}
In this section we have gathered the notations we use and the basic definitions we need in the
subsequent sections, see also \cite{Ma}.
\subsection{Notion of Graphs}
A directed weighted graph is a triple $G :=(V,\vec{E},b)$, where $V$ is the countable set of the vertices, $\vec{E}$ is the set of directed edges and $b : V \times V\to\left[0,\infty\right)  $ is a weight satisfying the following conditions:
\begin{itemize}
  \item $b(x,x)=0$ for all $x\in V$~~~(no loops)
  \item $b(x,y)>0$ iff $(x,y)\in \vec{E}$
\end{itemize}

In addition, we consider a measure on V given by a nonnegative real function
$$m:V\to \left]  0,\infty\right) .$$
The weighted graph is \textit{symmetric} if  for all $x,y\in V$, $b(x,y)=b(y,x)$,
  as a consequence $(x,y)\in \vec{E}\Rightarrow (y,x)\in\vec{E}$.\\
The graph is called  \textit{simple} if the weights $m$ and $b$ are constant equal to $1$ on $V$ and $\vec{E}$ respectively.\\

On a non symmetric graph we have several notions of {\it connexity}.

We fix the following notations:
\begin{itemize}
\item The set of {\it undirected edges} is defined by
$$E=\left\{\{x,y\},~(x,y)\in\vec{E}\text{ or } (y,x)\in\vec{E}\right\}.$$
\item for $x\in V,\,V_x^+ = \{y \in V ; (x, y) \in \vec{E}\}$
\item for $x\in V,\,V_x^- = \{y \in V ; (y,x) \in \vec{E}\}$
\item for $x\in V,\,V_x = V_x^+ \cup V_x^-=\{y \in V ; \{x,y\} \in E\}$.
\end{itemize}
\begin{defini}
 The degree of a vertex $x$ is denoted by $deg(x)$ and defined by:
 $$ deg(x)=\#V_x.$$
\begin{itemize}
\item  A {\it chain} from the vertex $x$ to the vertex $y$ in $G$ is a finite set of undirected edges $\{x_1,y_1\};~\{x_2,y_2\};..;\{x_n,y_n\},~n\geq 1$
$$x_1=x,~y_n=y \text{ and } x_i=y_{i-1}~~\forall~2\leq i\leq n.$$

\item A {\it path} between two vertices $x$ and $y$ in $V$ is a finite set of directed edges $(x_1,y_1);~(x_2,y_2);..;(x_n,y_n),~n\geq 1$ such that
$$x_1=x,~y_n=y \text{ and } x_i=y_{i-1}~~\forall~2\leq i\leq n$$
\item $G$ is called {\it weakly connected} if two vertices are always related by a chain.
\item $G$ is called {\it connected} if two vertices are always related by a path.
\item $G$ is called {\it strongly connected} if there is for all vertices $x,y$ a path from $x$ to $y$ and one from $y$ to $x$.
 \end{itemize}
\end{defini}
We assume in the following that the graph under consideration is weakly connected, locally finite
and satisfy:
$$\text{ for all } x\in V,~~\sum_{y\in V}b(x,y)>0.$$
\subsection{Functional spaces}
Let us introduce the following function spaces associated to the graph $G$.\\
The space of functions on the graph $G$ is considered as the space of complex
functions on V and is denoted by
$$\mathcal{C}(V)=\{f:V\to \mathbb{C} \}.$$
We denote by  $\mathcal{C}_c(V)$  its subset of finite supported functions. We consider for a measure $m$,
the space
$$\ell^2(V,m)=\{f\in \mathcal{C}(V), ~~\sum_{x\in V}m(x)|f(x)|^2<\infty\},$$
which is a Hilbert space when equipped by the scalar product given by
$$\< f,g\>=\sum_{x\in V}m(x)f(x)\overline{g(x)}.$$
The associated norm is given by:
$$\|f\|=\sqrt{\< f,f\>}.$$
\subsection{Laplacians and Kirchhoff's Assumption}
For a locally finite connected graph without loops, we introduce the weighted Laplacian $\Delta$
defined on $\mathcal{C}_c(V)$ by:
$$\Delta f(x)=\frac{1}{m(x)}\sum_{y\in V}b(x,y)\left( f(x)-f(y)\right) .$$

\textbf{Kirchhoff's Assumption $(\beta)$}: This assumption says that at each vertex the incoming
conductance equals the  outcoming conductance. If the Kirchhoff's Assumption is satisfied, the non
symmetric operator $\Delta$ enjoys parts of the self-adjoint theory. Defining for
$x\in V$, $\displaystyle\beta^-(x)=\sum_{y\in V}b(y,x)$ and $\displaystyle\beta^+(x)=\sum_{y\in V}b(x,y)$,
we will suppose in the sequel of this work that \\

\noindent
\begin{tabular}{lc}
($\beta$)&\hspace{4cm}$\forall x \in V,\;\beta^+(x)=\beta^-(x)$.
\end{tabular}\\

With this assumption, the formal adjoint of $\Delta$ has a simple expression:
 \begin{prop}
The formal adjoint $\Delta'$ of the operator $\Delta$ is defined on $\mathcal{C}_c(V)$ by:
$$\Delta' f(x)=\frac{1}{m(x)}\sum_{y\in V}b(y,x)\big(f(x)-f(y)\big).$$
\end{prop}
 In this situation, we have established, see \cite{Ma}, an explicit Green formula associated to the
 non symmetric Laplacian $\Delta$.
\begin{prop}(Green Formula)\label{grn}
Let $f$ and $g$ be two functions of $\mathcal{C}_c(V)$. They satisfy
$$\<\Delta f,g\>+\<\Delta' g,f\>=\sum_{(x,y)\in \vec{E}}b(x,y)\big( f(x)-f(y)\big) \big( \overline{g(x)-g(y)}\big).$$
\end{prop}

\section{m-accretiveness of the Laplacian}\label{2}
\subsection{First properties}
The  Hilbert  space theory  of  accretive  operators was  motivated  by  the semi-group theory and
the Cauchy  problem  for  systems  of  hyperbolic  partial differential  equations. It is an important
property for operators which cannot be studied in the framework of selfadjointness.
We establish conditions for the m-accretiveness (maximal accretiveness, see \dref{m-acc}) of
$\Delta$.
\begin{defi}
  The \emph{numerical range} of an operator $A$ with domain $D(A)$, denoted by  $W(A)$ is the non-empty
  subset of $\C$ defined by
$$W(A)=\{\<A f,f\>,~~f\in D(A),~\parallel f\parallel=1 \}.$$
\end{defi}

\begin{defi}
  Let $\mathcal{H}$ be a Hilbert space, an operator $A: D(A)\to H$ is said to be \emph{accretive} if for
  each $f\in D(A)$,
$$\Re\<A f,f\>\geq 0.$$
\end{defi}
\begin{lemm}\label{den}\cite[Lem.1.47]{Ou}
  Let $A$ be a densely defined accretive operator on $\mathcal{H}$. Then $A$ is closable, its closure
  $\overline{A}$ is accretive, and for every $\lambda > 0$, the range $\mathrm{Im} (\lambda+A)$ is
  dense in $\mathrm{Im} (\lambda+\overline{A})$.
\end{lemm}
\begin{prop}\label{imc}
Let $A:D(A)\to \mathcal{H}$ be a closed, densely defined and accretive operator, then
$$ \|(A+\lambda) f\| \geq \Re(\lambda) \|f\|,~~\forall f\in D(A),~\Re(\lambda)> 0$$
and $\mathrm{Im} (A+\lambda)$ is closed.
\end{prop}
\begin{demo}
Let $\Re(\lambda)> 0$, $\forall f\in D(A)$ we have
\begin{align*}
\|(A+\lambda) f\|\|f\|\geq |\<(A+\lambda )f,f\>|\geq &\Re\<(A+\lambda )f,f\> \\
\geq & \Re\<Af,f\>+\Re\<\lambda f,f\>\\
\geq & \Re(\lambda) \|f\|^2.
\end{align*}
Hence $ \forall f\in D(A),\,\|(A+\lambda) f\|\geq \Re(\lambda)\|f\|$ and Im$(A+\lambda)$ is closed.
\end{demo}
\begin{prop}In our situation, the Laplacian $\Delta$ with domain $D(\Delta)=\mathcal{C}_c(V)$
  is accretive, closable and for any scalar $\Re(\lambda) > 0$, $\mathrm{Im} (\Delta+\lambda)$ is dense in
  $\mathrm{Im} (\overline{\Delta}+\lambda)$.
\end{prop}
\begin{demo}
From the Green's formula, \pref{grn}, we have for any $f\in \mathcal{C}_c(V)$
\begin{align*}
\Re\<\Delta f,f\>=&\dfrac{1}{2}\left( \<\Delta f,f\>+\overline{\<\Delta f,f\>}\right) \\
=& \dfrac{1}{2}\left( \<\Delta f,f\>+\<\Delta' f,f\>\right)\\
=&\dfrac{1}{2}\sum_{(x,y)\in \vec{E}}b(x,y)\big| f(x)-f(y)\big|^2\\
\geq & 0.
\end{align*}
Therefore $\Delta$ is accretive. We deduce from \lref{den} that $\Delta$ is closable and
$\mathrm{Im}(\lambda+\overline{\Delta})=\overline{\mathrm{Im}(\lambda+\Delta)}$.
\end{demo}

We introduce the following notations (already introduced in \cite{A}): let
\begin{equation}\label{HB}
  H=\frac{1}{2}(\Delta+\Delta') \quad  B=\frac{1}{2}(\Delta-\Delta')
\end{equation}
be the symmetric and the skewsymmetric parts of $\Delta$, acting on the space of functions
with finite support.

Then, thanks to the assumption $(\beta)$, the operator $H$ is the Laplacian
on the symmetric graph with an edge $\{x,y\}$ weighted by the symmetric weight defined by
\begin{equation}\label{b'}
b'(x,y)=\displaystyle\frac{b(x,y)+b(y,x)}{2}~ \hbox{ for all } {x,y} \in E.
\end{equation}
\begin{defi}\label{m-acc}
  An accretive operator $A:D(A)\to \mathcal{H}$ is said to be \emph{m-accretive} if the left open half-plane
  is contained in the resolvent set $\rho(A)$ and we have for
  $\Re(\lambda)>0$, $$||(A+\lambda)^{-1}||\leq\dfrac{1}{\Re(\lambda)}.$$
\end{defi}

An m-accretive operator $A$ is maximal accretive, in the sense that $A$ is accretive and has no proper
accretive extension, \cite{Kts}. \\

In the following we give sufficient conditions for $\Delta$ to be m-accretive, based on previous works
applied to the real part of $\Delta$.
First, we recall a relating result, in the symmetric case,  essential selfadjointness to graphs with
constant weights on $V$.
\begin{theo}[Theorem 3.1 of \cite{To10}]\label{const}
  Let $(G,m,b')$ be an infinite weighted graph with a constant weight $m$ on $V$.  Then the Laplacian $H$
  is essentially selfadjoint.
\end{theo}

From the definition the adjoint operator $\Delta^*$, we can deduce:
$$D(\Delta^*)=\{f\in\ell^2(V,m),~~\Delta'f\in\ell^2(V,m)\}. $$
Using an idea in the proof of Theorem 3.1 of \cite{To10} , we prove the following Proposition for the non symmetric Laplacian.
\begin{prop}\label{constnonsym}Let $(G,m,b)$ be an infinite weighted graph with the constant weight $m$ on $V$.
  Then the Laplacian $\overline{\Delta}$ is m-accretive.
\end{prop}
\begin{demo}
We show that $\Delta^*+ 1$ is injective:
 Let $g\in \ell^2(V,m)$ be a function satisfying
$$\Delta'g+g=0, \text{ then } g=0:$$
Let us assume that there is a vertex $x_0$ such that $g (x_0 )>0$. The equality $\Delta'g+g=0$ implies that there exists at least one neighbouring vertex $x_1$ for which $g(x_0)<g (x_1 )$. We repeat the
procedure with $x_1$ ... Hence we build a strictly increasing sequence of
strictly positive real numbers $(g (x_n ))_n$. We deduce that the function
g is not in $\ell^2(V,m)$.
\end{demo}
\begin{rem}
If $\Delta$ is symmetric, $\overline{\Delta}$ is m-accretive if and only if $\Delta$ is essentially self-adjoint.
\end{rem}
The property of essential self-adjointness was extensively studied in the symmetric case and many tools related
to completeness were introduced to assure this property. In \cite{To10}, one of us proved that essential
self-adjointness followed from completeness for a certain metric of the graph with bounded degree.
In \cite{HKMW} the condition is related on completeness for an {\it intrinsic metric}. In \cite{AT}
we introduced the notion of $\chi$-completeness.

\subsection{$\chi$-completeness}
We have introduced this notion in \cite{AT} in the symmetric case. It assures the Laplacian
(and even the Gau\ss-Bonnet operator) to be essentially selfadjoint.
We suppose in this section that the graph is $\chi$-complete for the symmetric
conductance $b'$ defined in \eref{b'}. It means that there exists an increasing sequence of finite sets
$(B_n)_{n\in \N}$ such that
$V=\ccup\, B_n$ and there exist related functions $\chi_n$ satisfying the
following three conditions:
\begin{enumerate}
\item[  (i)] $\chi_n\in \mathcal C_c(V),\, 0\leq\chi_n\leq 1$
\item[ (ii)] $v\in B_n ~\Rightarrow ~\chi_n(v)=1$
\item[(iii)] $\displaystyle\exists C>0, \forall n\in \N,\, x\in V,\,
\frac{1}{m(x)}\sum_{y\in V_x }b'(x,y)|\chi_n(x)-\chi_n(y)|^2\leq C.$

\end{enumerate}

\begin{theo}\label{chi}Suppose that the graph $G=(V,m,b)$ is $\chi$-complete for the symmetric
  conductance $b'$, and that the asymmetry is controled in the following way
  \begin{equation}\label{thm-chi}
 \exists C >0,\, \forall x\in V,\,
  \frac{1}{m(x)}\sum_{y\in V_x}\frac{|b(x,y)-b(y,x)|^2}{b'(x,y)}\leq C
  \end{equation}
  then the nonsymmetric Laplacian $\overline{\Delta}$ is m-accretive.
\end{theo}

\begin{demo} We can suppose that the constants for $\chi$-completness and for
  \eqref{thm-chi} are the same. By \lref{den} \, $\overline{\Delta}$ is accretive
  and by \pref{imc} \,, its range is closed.
Suppose that $\overline{\Delta}$ is not m-accretive, it means that the deficiency of
  $\overline{\Delta}$, which is constant on the left halfplane, is not 0. For instance at $-1$ it gives
  $$\exists v\in \ell^2(V,m) , \forall x\in V,\, (\Delta'+1)v(x)=0.
  $$
  We remark that the operator $\Delta'$ is real so we can suppose that \emph{$v$ is real}.
  Let $\chi\in \mathcal C_c(V)$, such that  $0\leq \chi \leq 1$,  and calculate $\<\chi v,(\Delta'+1)(\chi v)\>$.\\
  First we remark that, as $\chi v\in \mathcal C_c(V)$ and has real values:
  \begin{multline*}
\<\chi v,(\Delta'+1)(\chi v)\>=\<(\Delta+1)(\chi v),\chi v\>=\<\chi v,(\Delta+1)(\chi v)\>=\\
\<\chi v,(H+1)(\chi v)\>\geq \|\chi v\|^2
  \end{multline*}
  On the other hand, using the equation satisfied by $v$ we have:
  \begin{multline*}
    (\Delta'+1)(\chi v)(x) =\\
\frac{1}{m(x)}\sum_{y\in V}b(y,x)\big(\chi(x)v(x)-\chi(y)v(y)\big)-
    \chi(x)\frac{1}{m(x)}\sum_{y\in V}b(y,x)\big(v(x)-v(y)\big)\\
    =\frac{1}{m(x)}\sum_{y\in V}b(y,x)v(y)\big(\chi(x)-\chi(y)\big)
  \end{multline*}
  it gives
  \begin{multline*}
    \<\chi v,(\Delta'+1)(\chi v)\>=\sum_{x\in V}\chi(x)v(x)\sum_{y\in V}b(y,x)v(y)\big(\chi(x)-\chi(y)\big)\\
    =\frac{1}{2}\sum_{x\in V,y\in V}v(x)v(y)\big(\chi(x)-\chi(y)\big)\big(b(y,x)\chi(x)-b(x,y)\chi(y)\big)\\
    =\frac{1}{2}\sum_{x\in V,y\in V}v(x)v(y)\big(b(y,x)\chi(x)^2+b(x,y)\chi(y)^2-\chi(x)\chi(y)(b(x,y)+b(y,x))\big)\\
    =\frac{1}{2}\sum_{x\in V,y\in V}v(x)v(y)\big(b(y,x)\chi(x)^2+b(x,y)\chi(y)^2-2\chi(x)\chi(y)b'(x,y)\big).
  \end{multline*}
 We use then that $2|v(x)v(y)|\leq v(x)^2+v(y)^2$, it gives
   \begin{multline*}
     \<\chi v,(\Delta'+1)(\chi v)\>\leq \\
     \frac{1}{2}\sum_{x\in V}v(x)^2
     \sum_{y\in V}|b(y,x)\chi(x)^2+b(x,y)\chi(y)^2-2\chi(x)\chi(y)b'(x,y)|.
   \end{multline*}
   We see that
   $$ \hbox{if } \chi(x)=\chi(y)=1 , \hbox{ then } \big(b(y,x)\chi(x)^2+b(x,y)\chi(y)^2-2\chi(x)\chi(y)b'(x,y)\big)=0. $$
  Moreover
   \begin{multline*}
     b(y,x)\chi(x)^2+b(x,y)\chi(y)^2-2\chi(x)\chi(y)b'(x,y)=\\
   b'(x,y)\big(\chi(x)-\chi(y)\big)^2+\frac{b(y,x)-b(x,y)}{2}\big(\chi(x)^2-\chi(y)^2\big)
   \end{multline*}
   We remark that 
   \begin{multline*}
   |\chi(x)^2-\chi(y)^2|=|\chi(x)-\chi(y)|.\big(\chi(x)+\chi(y)\big)\leq 2|\chi(x)-\chi(y)|,
   \end{multline*}
 which implies
   \begin{multline*}
   \sum_{y\in V}|\frac{b(y,x)-b(x,y)}{2}\big(\chi(x)^2-\chi(y)^2\big)|
   \leq\\
   \sqrt{\sum_{y\in V}\frac{|b(y,x)-b(x,y)|^2}{b'(x,y)}}
   \sqrt{\sum_{y\in V}b'(x,y)|\chi(x)-\chi(y)|^2}
   \end{multline*}
   Applying this calculation to $\chi=\chi_n$ we have then, because of the hypothesis \eref{thm-chi}
   \begin{equation*}
  \sum_{y\in V}|b(y,x)\chi_n(x)^2+b(x,y)\chi_n(y)^2-2\chi_n(x)\chi_n(y)b'(x,y)|\leq 2Cm(x)
   \end{equation*}
   and finally
   \begin{equation*}
     \|\chi_n v\|^2\leq \<\chi_n v,(\Delta'+1)(\chi_n v)\>\leq 2C\sum_{x\in W_n}m(x) v(x)^2
   \end{equation*}
   where $W_n=V\setminus V_n$ and $V_n=\{x\in B_n,\, \forall y\in V_x, \chi_n(y)=1\}$.
   Thus $\lim_{n\to\infty}\sum_{x\in W_n}m(x) v(x)^2=0$ (as $v\in \ell^2(V,m)$), we conclude
   that $v=0$.

   Thus the deficiency of $\overline{\Delta}$ is 0 on the left halfplane, and
   we conclude by \pref{imc}.
\end{demo}
\begin{rem}
  The results of essential self-adjointness for the Laplacian in the symmetric case
  give results for $H$ in our case. Indeed the hypothesis of $\chi$-completeness for
  the weight $b'$ gives that $H$ is essentially selfadjoint.
\end{rem}
\begin{exem}\label{tree}
 Let us consider an infinite simple tree $T$ with increasing degree, see Fig.1, we suppose that
\begin{equation}\label{hyp-dir}
\forall x\in V \quad \#V_x^+\setminus (V_x^+ \cap V_x^-)=\#V_x^-\setminus (V_x^+ \cap V_x^-)=1.
\end{equation}

\begin{figure}[ht]
\includegraphics*[height=10cm,width=14cm]{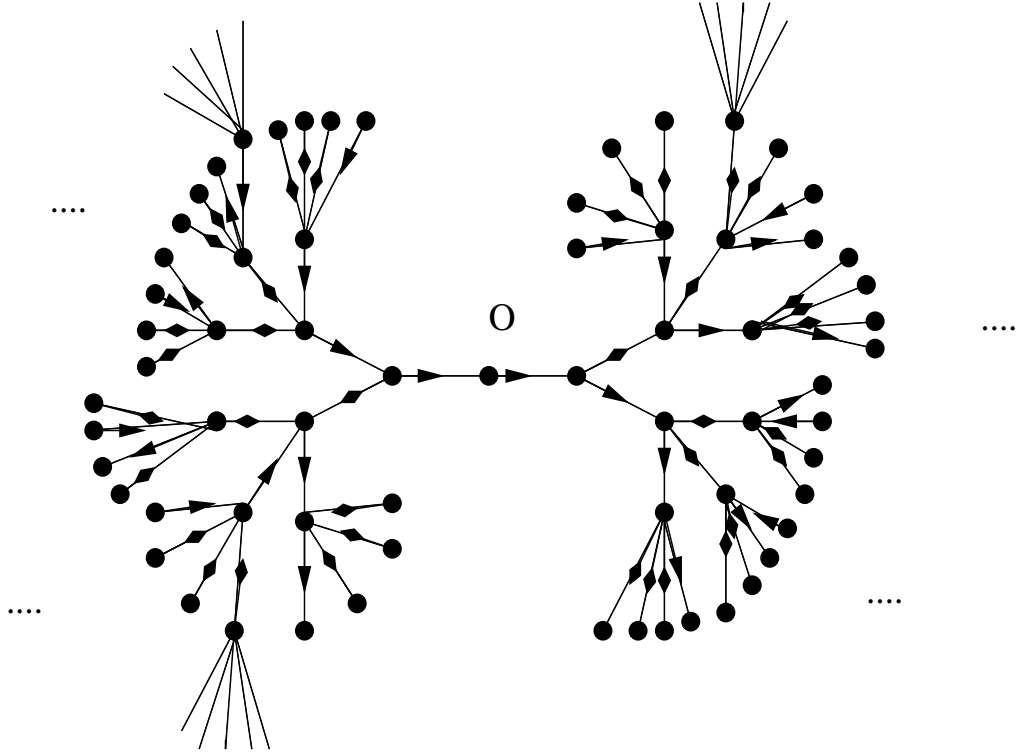}
\caption{a tree with increasing degree}
\end{figure}

It is shown in \cite[Example 9]{AT} that $T$, considered as a symmetric graph,
is $\chi$-complete. The hypothesis \eqref{hyp-dir} assures clearly the Kirchhoff
law $(\beta)$, we can see also that the property \eqref{thm-chi} is satisfied, in fact
we have for all $x\in V$
\begin{align*}
  \frac{1}{m(x)}\sum_{y\in V_x}\frac{|b(x,y)-b(y,x)|^2}{b'(x,y)}&=
 \sum_{y\in V_x^+\setminus (V_x^+ \cap V_x^-)}
  +\sum_{y\in V_x^-\setminus (V_x^+ \cap V_x^-)}\\
  &=\#V_x^+\setminus (V_x^+ \cap V_x^-)+\#V_x^-\setminus (V_x^+ \cap V_x^-)\\
  &=2.
\end{align*}

 More generally, we can suppose that
 $$\#V_x^+\setminus (V_x^+ \cap V_x^-)=\#V_x^-\setminus (V_x^+ \cap V_x^-)\leq M.$$
 But, using the fact that the degree is not bounded, we can also construct
 a graph with $\#V_x^+\setminus (V_x^+ \cap V_x^-)=\#V_x^-\setminus (V_x^+ \cap V_x^-)$
 not bounded. This gives a $\chi$- complete graph which does not satisfy the property
 \eqref{thm-chi}.
\end{exem}
\begin{rem}
  In \cite{syl} the authors give different criteria for $\chi$-completness on weighted graphs.
  Let us consider a $\chi$-complete graph (for the symmetrized weight $b'$) following
  Proposition 5.7 or Theorem 5.11 of \cite{syl}.
  To obtain the m-accretiveness of $\overline{\Delta}$, it is then sufficient that the non
  symmetric graph satisfies moreover the Kirchhoff Assumption $(\beta)$ and the property \eqref{thm-chi}.
  This is assured if we suppose for instance that
$\forall x\in V,~y\in V_x^+\cap V_x^-,~b(x,y)=b(y,x)$ and
\[\exists M>0,\,\forall x\in V~,~~~
 \sum_{y\in V_x^+\setminus (V_x^+ \cap V_x^-)}b(x,y)
  =\sum_{y\in V_x^-\setminus (V_x^+ \cap V_x^-)}b(y,x)\leq M.\]
\end{rem}
\begin{exem}\label{gg}
  Let us consider the following infinite weighted graph $G$, see Fig.2, with (almost) constant degree.
  We denote the origin by $x_0$ and  by $S_n$ the spheres for the combinatoric distance of the symmetric
  underlying graph:
\[d_{\rm comb}(x_0,x)=\inf\{k;\,\exists \gamma=(x_0,\dots,x_k) \hbox{ a chain such that }x_k=x\}\]
So $S_n=\{x\in V,~ d_{\rm comb}(x_0,x)=n\}=\{x_n,y_n\}$.\\

 \begin{figure}[ht]
\includegraphics*[height=4cm,width=10cm]{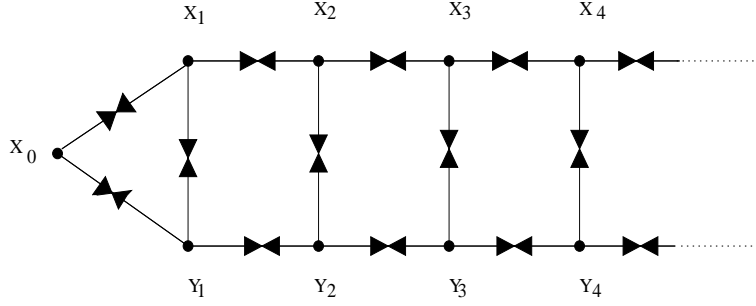}
\caption{a graph with almost constant degree}
 \end{figure}

To define the weights, we take $k\geq 0$ and fix
\begin{itemize}
\item $b(x_0,x_1)=b(y_1,x_0)=k+2$ and $b(x_0,y_1)=b(x_1,x_0)=k$ and for $n\geq 1$
\item $b(x_n,x_{n+1})=(n+1)^2+(n+1)$ and $b(x_{n+1},x_n)=(n+1)^2-(n+1)$
\item $b(y_n,y_{n+1})=(n+1)^2-(n+1)$ and $b(y_{n+1},y_n)=(n+1)^2+(n+1)$
\item $b(x_n,y_n)=n-1$ and $b(y_n,x_n)=n+1$
\item $m(x_0)=1$, $m(x_n)=m(y_n)=\sqrt{n}$.
\end{itemize}
We can see that the property \eqref{thm-chi} is satisfied, in fact for all
$x\in S_n=\{x_n,y_n\}$
\begin{align*}
  \frac{1}{m(x)}\sum_{y\in V_x}\frac{|b(x,y)-b(y,x)|^2}{b'(x,y)}&= \frac{1}{\sqrt{n}}\left(
 \sum_{y\in S_{n+1}\cap V_x}4
  +\sum_{y\in S_{n-1} \cap V_x}4+\sum_{y\in S_{n}\cap V_x}\dfrac{2^2}{n}\right) \\
  &\leq \dfrac{12}{\sqrt{n}}\leq 12.
\end{align*}
We show now that $G$ is $\chi$-complete from the criterion given in Theorem 5.11 of \cite{syl}:
we remark that the set $S_n^+$ and $S_n^-$ introduced in \cite{syl} coincide with $S_n$ and for $x\in S_n$ the weighted
degree is constant and for $n\geq 2$:
\begin{gather*}
  a_n^+=\sup_{x\in S_n }\frac{1}{\sqrt{n}}\left(\sum_{y\in S_{n+1}}b'(x,y)\right)=\dfrac{(n+1)^2}{\sqrt{n}}\hbox{ , }
  a_n^-=\sup_{x\in S_n }\frac{1}{\sqrt{n}}\left(\sum_{y\in S_{n-1}}b'(x,y)\right)=n^\frac{3}{2}\\
\Rightarrow
  \sum_{n=0}^\infty\dfrac{1}{\sqrt{a_n^++a_{n+1}^-}}\geq\sum_{n=2}^\infty\dfrac{1}{\sqrt{2(n+1)^\frac{3}{2}}}=\infty.
\end{gather*}
Thus, this graph satisfies the hypothesis of \tref{chi}.
\end{exem}
\section{Relations between $\Delta$ and $H$}We study here the relation between the two hypothesis:
m-accretiveness for $\overline{\Delta}$ and essential selfadjointness for $H$.
\subsection{From $H$ to $\Delta$}
\begin{theo}Let $(G,m,b)$ be an infinite weighted graph and $\Delta=H+B$ the decomposition
  of the combinatorial Laplacian of $G$ decomposed in symmetric and skewsymmetric part as
  in \eref{HB}. Then if $H$ is essentially selfadjoint and if $B$ is bounded,
  then $\overline{\Delta}$ is m-accretive.
\end{theo}
\begin{demo}As $B$ is bounded, $\overline\Delta$ and $\overline H$ have the same domain and
  \[\hbox{on } \mathcal D(\overline\Delta)=\mathcal D(\overline H) :\;\overline\Delta=\overline H +B.
  \]
  As $H$ is non-negative and essentially selfadjoint, for any scalar $\lambda,\Re(\lambda)>0$,
  the operator $\overline H+\lambda$ is invertible on $\mathcal D(\overline H)$ and
  \[\|(\overline H+\lambda)^{-1}\|\leq\frac{1}{\Re(\lambda)}.
  \]
  Now, let $\lambda,\Re(\lambda)>0$, we write
  \[\overline\Delta+\lambda=\overline H+\lambda+B=(I+B(\overline H+\lambda)^{-1})(\overline H+\lambda)
  \]
  But $\|B(\overline H+\lambda)^{-1})\|\leq \frac{\|B\|}{\Re(\lambda)}$ so
  $  \Re(\lambda)>\|B\|\Rightarrow (\overline\Delta+\lambda)$ invertible and
  \[ \|(\overline\Delta+\lambda)^{-1}\|
\leq\frac{1}{(1-\|B\|/\Re(\lambda))\Re(\lambda)}=\frac{1}{\Re(\lambda)-\|B\|}.
\]
So we have that the set of $-\lambda,\Re(\lambda)>\|B\|$ is included in
the resolvent set of $\overline\Delta$ but, on the other hand $\overline\Delta$
is accretive which implies, by Theorem V.3.2 of \cite[p.268]{Kts}, that its deficiency is constant
on the set of $\lambda,\Re(\lambda)<0$, as a conclusion this deficiency is zero on
this set.
Finally, for all $\lambda,\Re(\lambda)>0$ and $f\in C_c(V)$, because the real part of $\Delta$ is non-negative
\[
\Re(\lambda) <f,f>\leq |<(\lambda+\Delta)f,f>|\leq \|(\lambda+\Delta)f\|.\|f\|
\]
and this is also true on $\mathcal D(\overline \Delta)$, then, as we already know that $(\lambda+\overline\Delta)$ is invertible,
\[\|(\lambda+\overline\Delta)^{-1}\|\leq\frac{1}{\Re(\lambda)}.
\]
\end{demo}
\begin{rem}\label{gamma}We have introduced in \cite{A} an hypothesis $(\gamma)$ that assures $B$
  to be bounded (and $\Delta$ to be sectorial), namely

\textbf{Assumption $(\gamma)$}:\\
\begin{equation*}
\exists~M>0,~\forall~x\in V,~\sum_{y\in V}\mid b(x,y)-b(y,x)\mid \leq M m(x)
\end{equation*}
We see easily that if the assumption $(\gamma)$ is satisfied then \eref{thm-chi} is also satisfied:
 as the weight $b$ is non-negative, we have always
 \begin{align*}
   |b(x,y)-b(y,x)|&\leq (b(x,y)+b(y,x))=2b'(x,y)
 \end{align*}
 and thus
 \begin{align*}
 \sum_{y\in V_x}\frac{|b(x,y)-b(y,x)|^2}{b'(x,y)}&\leq 2\sum_{y\in V_x}|b(x,y)-b(y,x)|\\
&\leq 2M m(x).
 \end{align*}

\end{rem}
\begin{exem}
  The graph considered in the  example \ref{gg} satisfies the property \eqref{thm-chi} and does not satisfies
  $(\gamma)$. In fact we have for all $x\in S_n$\\
\begin{align*}
\displaystyle{\sum_{y\in V_x}\mid b(x,y)-b(y,x)\mid} &= \displaystyle{\sum_{y\in S_n}\mid b(x,y)-b(y,x)\mid +\sum_{y\in S_{n+1}}\mid b(x,y)-b(y,x)\mid} \\
&+\displaystyle{\sum_{y\in S_{n-1}}\mid b(x,y)-b(y,x)\mid}\\
&=2+ 2(n+1)+2n
\end{align*}
which can not be controled by $m(x)=\sqrt n$.
\end{exem}
\begin{rem}\label{BHborne}
  The last theorem can be extended in a situation more general than sectoriality (see the definition
  in  \sref{sector}), namely when $B$ is bounded in $H$-norm with a relative norm sufficiently
  small. More precisely we suppose that there exist two constants $C>0$ and $0<a<1/2$ such that
  \[ \forall f\in C_c(V)\quad\|B(f)\|\leq C\|f\|+a\|H(f)\|
  \]
  Then for a real $\lambda>0$ we have
  \begin{multline*}
\forall f\in C_c(V)\; \|B(\overline H+\lambda)^{-1}f\|
\leq (C+a\lambda)\|(\overline H+\lambda)^{-1}f\|+a\|f\|\\
\Rightarrow \|B(\overline H+\lambda)^{-1}\|\leq \frac{C+a\lambda}{\lambda}+a\leq \frac{C}{\lambda}+2a
  \end{multline*}
  This can be made smaller than 1 for $\lambda$ large enough and then, by the same argument,
  the deficiency of $\overline\Delta$ must be zero on all the left halfspace (notice that under these hypothesis $B$ is also $\Delta$-bounded, with relative norm $\frac{a}{1-a}$).
\end{rem}
\begin{theo}\label{Hb}Let $(G,m,b)$ be an infinite weighted graph and $\Delta=H+B$ the decomposition
  of the combinatorial Laplacian of $G$ decomposed in symmetric and skewsymmetric part as
  in \eref{HB}. If $H$ is essentially selfadjoint and if $B$ is relatively bounded with respect
  to $H$ with relative norm smaller than $1/2$, then $\overline{\Delta}$ is m-accretive.
\end{theo}
\begin{rem}The hypothesis \eref{thm-chi} gives that $B$ is relatively bounded with respect to $H$.
  Indeed for any $f\in\mathcal C_c(V)$
  \begin{multline*}
    \|B(f)\|^2=\sum_{x\in V} \dfrac{1}{m(x)}\lvert\sum_{y\in V}\frac{b(x,y)-b(y,x)}{2}(f(x)-f(y))\rvert^2\\
    \leq\sum_{x\in V}\dfrac{1}{m(x)}\sum_{y\in V}\frac{\lvert b(x,y)-b(y,x)\rvert^2}{4b'(x,y)}
    \sum_{y\in V}b'(x,y)\lvert f(x)-f(y)\rvert^2\\
    \leq\sum_{x\in V}\frac{C}{4} \sum_{y\in V}b'(x,y)\lvert f(x)-f(y)\rvert^2 =\frac{C}{2}\<H(f),f\>\\
    \leq \Big(\frac{C^2}{4}\|f\|^2+\frac{1}{4}\|H(f)\|^2\Big).
  \end{multline*}
  Thus \tref{chi} is a corollary of \tref{Hb}.
\end{rem}
\begin{rem}If $\overline{\Delta}$ is m-accretive then, by definition, the set of $\lambda,\Re(\lambda)<0$
  is included in the resolvent set of $\overline\Delta$, we have thus
  \[\overline{\Delta} \hbox{ m-accretive}\Rightarrow \sigma(\overline\Delta)\subset
  \{\lambda\in \C,\,\Re(\lambda)\geq 0\}
  \]
\end{rem}
We can study if the hypothesis ``$H$ is essentially selfadjoint'' is necessary.\\

\subsection{From $\Delta$ to $H$}
\begin{prop}Let $(G,m,b)$ be an infinite weighted graph and $\Delta=H+B$ the decomposition
  of the combinatorial Laplacian of $G$ decomposed in symmetric and skewsymmetric part as
  in \eref{HB}. If $B$ is bounded and $\overline{
    \Delta}$ is m-accretive, then $\overline{\Delta'}$ is m-accretive
  and $H$ is essentially selfadjoint.
\end{prop}
\begin{demo}On $C_c(V)$ we have that $\Delta'=\Delta-2B$. The operator $\Delta'$ is
  accretive, as $B$ is bounded
  $\mathcal D(\overline\Delta)=\mathcal D(\overline\Delta')$ and for any $\lambda,\Re(\lambda)> 0$
  \[\Delta'+\lambda=\Delta+\lambda-2B=(I-2B(\Delta+\lambda)^{-1})(\Delta+\lambda)
  \]
  then $\overline{\Delta'}+\lambda$ is invertible for $\Re(\lambda)$ large enough, so
  $\overline{\Delta'}$ is m-accretive.
  In the same way we have, as $H=\Delta-B$ so $\mathcal D(\overline H)=\mathcal D(\overline\Delta)$
  and
  \[H+\lambda=\Delta-B+\lambda=(I-B(\Delta+\lambda)^{-1})(\Delta+\lambda),
  \]
  that $(\overline H+\lambda)$ is invertible for $\lambda$ real large enough and $H$ is essentially
  selfadjoint.
\end{demo}
\begin{rem}In the same way as in \rref{BHborne} we can extend this result for
  $B$ bounded in $\Delta$-norm. We obtain that if $\overline{\Delta}$ is m-accretive and
  $B$ is $\Delta$-bounded with a relative norm strictly smaller than $1/2$ then
  $H$ is essentially selfadjoint.
\end{rem}
\section{Sectoriality}\label{sector}
In \cite{A} we have studied the sectoriality of $\Delta$, we generalize here these results.
It was McIntosh, see \cite{Mc}, who initiated and developed a theory of functional calculus
for a less restricted large class of operators, namely sectorial operators.

\begin{Def}
  Let $\mathcal{H}$ be a Hilbert space, an operator $A: D(A)\to \mathcal{H}$ is said to be sectorial
  if $W(A)$ lies in a sector
$$S_{a,\theta}:=\{z\in\mathbb{C}, ~|\Im(z)| \leq \tan \theta (\Re z -a)\}:=\{z\in\mathbb{C},~\mid arg(z-a)\mid \leq \theta \}$$
for some $a\in \mathbb{R}$, called vertex of $S_{a,\theta}$ , and $\theta\in \left[  0,\frac{\pi}{2}\right) $, called semi-angle of $S_{a,\theta}$ (thus $A-a$ is accretive).
The operator $A$ is said to be m-sectorial, if it is  sectorial and if $A-a$ is m-accretive.
\end{Def}

We have used in \cite{A} that under the assumption $(\gamma)$ (see \rref{gamma}) the Laplacian is
sectorial.
More generally, we have
\begin{pro}Let $(G,m,b)$ be an infinite weighted graph and $\Delta=H+B$ the decomposition
  of the combinatorial Laplacian of $G$ decomposed in symmetric and skewsymmetric part as
  in \eref{HB}. If the assymmetry of the weight $b$ satisfies the property \eref{thm-chi} then
  $\Delta$ is sectorial.
\end{pro}
\begin{demo}  Let $f\in \mathcal{C}_c(V)$ with $\|f\|=1$, using the Cauchy-Schwarz inequality we have
  \begin{align*}
  2\left|(Bf,f)\right|&=\left|\sum_{x\in V}f(x)\sum_{y\in V}(b(x,y)-b(y,x))(f(x)-f(y))\right|\\
  &= \left|\sum_{x\in V}f(x)\sum_{y\in V}\dfrac{b(x,y)-b(y,x)}{\sqrt{b'(x,y)}}\sqrt{b'(x,y)}(f(x)-f(y))\right|\\
  &\leq \sum_{x\in V}|f(x)|\left( \sum_{y\in V}\dfrac{\left|b(x,y)-b(y,x)\right|^2}{b'(x,y)}\right)^{\frac{1}{2}}\left(  \sum_{y\in V}b'(x,y)|f(x)-f(y)|^2\right)^\frac{1}{2}\\
  & \leq \sqrt C \left( \sum_{x\in V}m(x)|f(x)|^2\right)^{\frac{1}{2}}\left(\sum_{x\in V} \sum_{y\in V}b'(x,y)|f(x)-f(y)|^2\right)^\frac{1}{2}\\
  & \leq \sqrt C \|f\|(Hf,f)^{\frac{1}{2}}\\
  &\leq 1+\frac{C}{4}(Hf,f).
  \end{align*}
\end{demo}
\begin{pro}
Suppose that the graph $(G,m,b)$ is $\chi$-complete for the symmetric
  conductance $b'$, and that the assymmetry of the weight $b$ satisfies the property \eref{thm-chi} then the nonsymmetric Laplacian $\overline{\Delta}$ is m-sectorial.
\end{pro}

\section{The heat semigroup}
The property of m-accretivity can be used to generate strongly continuous semigroups.
We recall the theorem of Hille-Yosida. It gives, on Banach spaces, a complete
characterization of generators of semigroups with at most exponential growth,
we refer here to \cite{Re}.
\subsection{Existence of a heat semigroup}
\begin{theo}[Hille-Yosida]
  Let A be an operator in the Banach space $X$. Then $A$ is the infinitesimal generator of a $C_0$
  semigroup $T(t)$ satisfying $\|T (t)\| \leq M \exp(wt)$ if and only if the following two conditions
hold:
\begin{itemize}
\item $D(A)$ is dense and $A$ is closed.
\item Every real number $\lambda > w$ is in the resolvent set of A and
  $$\|(A-\lambda)^{-n}\|\leq \dfrac{M}{(\lambda-w)^{-n}}, \text{ for every }n\in \mathbb{N}.$$
\end{itemize}
\end{theo}
The assumptions of the Hille-Yosida Theorem are easier to achieve when $M=1$, then
the semigroup is said to be quasicontractive (and contractive if we can take $w=0$).
It is known as the Lumer-Phillips Theorem.
\begin{theo}[Lumer-Phillips]
  Let $A$ be a linear operator on a Hilbert space $\mathcal{H}$. If
  \begin{enumerate}
  \item $D(A)$ is dense
  \item $\Re(x,Ax)\leq w (x,x)$ for  $x\in D(A)$
  \item there exists $\lambda_0 > w$ such that $A-\lambda_0 I$ is onto.
  \end{enumerate}
  Then $A$ is the generator of a
strongly continuous one-parameter quasicontraction semigroup and $\|\exp(t A)\| \leq \exp(t w)$.
\end{theo}
  We can apply this to $A=-\overline{\Delta}$:
  \begin{theo}Let $(G,m,b)$ be an infinite weighted graph and $\overline{\Delta}$ its combinatorial
    Laplacian. If $\overline{\Delta}$ is m-accretive, then $-\overline{\Delta}$ is the generator of a
strongly continuous one-parameter contraction semigroup ({\it i.e.} $\|\exp(-t\overline{\Delta})\|\leq 1$).
  \end{theo}
  Using the results of \cite{Z} we have in the same way
\begin{theo}Let $(G,m,b)$ be an infinite weighted graph and $\overline{\Delta}$ its combinatorial
  Laplacian. If $\overline{\Delta}$ is m-sectorial with angle $\theta$ and vertex $a$, then
  $-\overline{\Delta}$ is the generator of an holomorphic semigroup on a sector with angle
  $\pi/2-\theta$ and vertex 0.
  \end{theo}
Moreover, on this situation one can apply the preceding result on $(\overline{\Delta}-a)$.
\subsection{Fast contractivity}
For a non symmetric graph $G$, we can estimate bounds on the real part of the numerical range of $\Delta$
in terms of a Cheeger constant. We restrict here in the case where the weight on vertices is constant equal
to 1 and consider the definition of the Cheeger constant given by Dodziuk, applied on the
symmetrized graph $(G,1,b')$.
\begin{defi}[\cite{Dod}]
Let us consider a weighted symmetric graph $(G,1,b')$, the Cheeger constant $h(V)$ is defined by
\[h(V)=\inf_{U\subset V \atop \mathrm{finite}}\displaystyle{\dfrac{\displaystyle{\sum_{x\in U,~y\in V\setminus U\atop y\in V_x}}\sqrt{b'(x,y)}}{\#U}}.
\]
\end{defi}
The following theorem is a consequence of Theorem 3.7 of \cite{M} and Theorem 3.1 of \cite{Dod}
(where deg is the combinatorial degree).
\begin{theo}
Suppose that $\sup_{x\in V} \mathrm{deg}(x)=M< \infty$. Then the real part of
$W(\Delta)$ satisfies
$$\inf_{z\in W(\Delta)}\Re(z)\geq\dfrac{h^2(V)}{2M}=\lambda_0.$$
\end{theo}
\begin{prop}\label{sect}
  Let $G$ be a graph with bounded degree, satisfying the property $\eqref{thm-chi}$ and $m=1$ on $V$.
  If $h(V)>0$, then $\Delta$ is m-sectorial with vertex $a\geq 0$.
\end{prop}
Indeed $\overline{\Delta}$ is m-accretive because of  \pref{constnonsym}. We remark that
$\overline{\Delta}-\lambda_0$ is also m-accretive.
\begin{exem}
  Consider the graph of Example \ref{gg}, but now with constant weight on vertices:  $m=1$.
  We define $V_n=\{x_0,x_1,y_1,...,x_n,y_n\}$  and remark that
 \begin{align*}
h(V)&=\inf_{n\in \mathbb{N}}\displaystyle{\dfrac{\displaystyle{\sum_{x\in \{x_n,y_n\}\atop y\in\{x_{n+1},y_{n+1}\}}}\sqrt{b'(x,y)}}{\#U}}\\
 &=\inf_{n\in \mathbb{N}}\dfrac{\sqrt{b'(x_n,x_{n+1})}+\sqrt{b'(y_n,y_{n+1})}}{\#U}\\
&=\inf_{n\in \mathbb{N}}\dfrac{2(n+1)}{2n+1}\\
&=1.
\end{align*}
Hence $$\inf_{z\in W(\Delta)}\Re(z)\geq \dfrac{1}{6}.$$
Then, applying \pref{sect}, there is $\theta\in(0,\pi/2)$ such that the numerical range of $\Delta$ lies in the
sector $\{z\in\mathbb{C},~\mid arg(z)\mid \leq \theta \}$ and $\overline{\Delta}$ is m-sectorial.

But we can say also that $(\overline{\Delta}-\dfrac{1}{6})$ is m-accretive, it gives
\[\|\exp(-t\overline{\Delta})\|\leq e^{-\frac{t}{6}}.
\]
\end{exem}

\textbf{Acknowledgments}:
The  author  Marwa  Balti  enjoyed  a  financial  support  from  the  Program  {\it D\'efiMaths} of
the  Federation  of  Mathematical Research of the "Pays de Loire" during her visits to the Laboratory of
Mathematics Jean Leray of Nantes (LMJL). Also, the three authors would like to thank the Laboratory of
Mathematics  Jean  Leray  of  Nantes  (LMJL)  and  the  research  unity  (UR/13ES47) of Faculty of
Sciences of Bizerta (University of Carthage) for their continuous financial support.

\end{document}